 \documentclass[12pt]{amsart}
\usepackage{amssymb}
\usepackage{epsf}
\usepackage{xy}

\numberwithin{equation}{section}

\theoremstyle{plain}
 \newtheorem{theo}{Theorem}[section]
 \newtheorem{prop}[theo]{Proposition}
 
 \newtheorem{coro}[theo]{Corollary}

\theoremstyle{definition}
 \newtheorem{defi}[theo]{Definition}

 \newtheorem{rem}[theo]{Remark}

 \newcommand{\gap}{\mathrm{gap}}       
\newcommand{\In}{\mathrm{In}}
           
\newcommand{\m}{\mathrm{Min}}
\newcommand{\M}{\mathrm{Max}}
\newcommand{\rac}{\mathrm{root}}
 
 \newcommand{\rto}{\rightarrow}

\newcommand{\calP}{\mathcal{P}}

 \newcommand{\calT}{\mathcal{T}}

\newcommand{\calS}{\mathcal S}

 \def\TT{{\mathbb T}}
 \def\YY{{\mathbb Y}}
 
\def\Lie{{ \mathcal {L}ie}} 
\def\As{{ \mathcal {A}s}} 
  \def\PL{{ \mathcal {PL}}}
\def\NAP{{\hbox{NAP}}}

\xyoption{all}
\CompileMatrices

\begin{document}

 \title[The non-symmetric operad pre-Lie is free]{The non-symmetric operad pre-Lie is free}

 \author[N.~Bergeron and M.~Livernet]{Nantel Bergeron
 and Muriel Livernet}

 \address[Nantel Bergeron]
 {Department of Mathematics and Statistics\\York University\\
 Toronto, Ontario M3J 1P3\\
Canada}
 \email{bergeron@mathstat.yorku.ca}
 \urladdr[Nantel Bergeron]{http://www.math.yorku.ca/bergeron}

 \address[Muriel Livernet]
 {Universit\'e Paris13\\CNRS, UMR 7539 LAGA \\ 
99,  Avenue
   Jean-Baptiste Cl\'ement\\ 93430 Villetaneuse\\France}
\email{livernet@math.univ-paris13.fr}
\urladdr[Muriel Livernet]{http://www.math.univ-paris13.fr/~livernet}

 \date{}

 \thanks{Livernet supported by the Clay Mathematical Institute and hosted by MIT}
 

 \begin{abstract}We prove that the pre-Lie operad is a free non-symmetric operad.
 \end{abstract}

 \maketitle

\vspace*{-0.3cm}\begin{figure}[!h]
\parbox{350pt}{\tiny\tableofcontents}
\end{figure}

 
 \section*{Introduction}

Operads are a specific tool for encoding type of algebras. For instance there are operads encoding associative algebras, commutative and associative algebras,
Lie algebras, pre-Lie algebras, dendriform algebras, Poisson algebras and so on. A usual way of describing a type of algebras is by giving the generating operations and 
the relations among them. For instance a Lie algebra $L$ is a vector space together with a bilinear product, the bracket (the generating operation) satisfying the relations
$[x,y]=-[y,x]$ and $[x,[y,z]]+[y,[z,x]]+[z,[x,y]]=0$ for all $x,y,z\in L$. The vector space of all operations one can perform on $n$ distinct variables
in a Lie algebra is $\Lie(n)$, the building block of the symmetric operad $\Lie$. 
Composition in the operad corresponds to composition of operations. The vector space
$\Lie(n)$ has a natural action of the symmetric group, so it is a symmetric operad. The case of associative algebras can be considered in two different ways. 
An associative algebra $A$ is a vector space together with a product satisfying the relation $(xy)z=x(yz)$. The vector space of all operations 
one can perform on $n$ distinct variables
in an associative algebra is $\As(n)$, the building block of the symmetric operad $\As$. The vector space $\As(n)$ has for basis 
the symmetric group $S_n$. But, in view of the relation, one can look also at  the vector space of all order-preserving operations 
one can perform on $n$ distinct ordered variables in an associative algebra: this is a vector space of dimension 1 generated by the only operation 
$x_1\cdots x_n$. So the non-symmetric operad $\widetilde{\As}$ describing associative algebras is 1-dimensional for each $n$: this is the terminal object in the 
category of non-symmetric operads. 

Here is the connection between symmetric and non-symmetric operads. A symmetric operad $\calP$ starts with a graded
vector space $(\calP(n))_{n\geq 0}$ together with an action of the symmetric group $S_n$ on $\calP(n)$ for each $n$. 
This data is called a symmetric sequence or an $\mathbb S$-module
or a vector species. There is a forgetful functor from the category of vector species to the category of graded vector spaces, forgetting the action of the symmetric group. This functor has
a left adjoint $\calS$ which corresponds to tensoring by the regular representation of the symmetric group. 
A symmetric (non-symmetric) operad is a monoid in the category of vector species (graded vector spaces). Again there is a forgetful functor from the category of symmetric operads to the category of non-symmetric operads admitting a left adjoint $\calS$. 
The symmetric operad $\As$ is the image of the non-symmteric operad  $\widetilde{\As}$ by $\calS$. It is clear that $\Lie$ is not in the image of $\calS$ since the Jacobi relation
does not respect the order of the variables $x<y<z$ nor the anti-symmetry relation. Still one can regard $\Lie$ as a non-symmetric operad applying 
the forgetful functor. Salvatore and Tauraso proved in \cite{SalTau08} that the operad $\Lie$ is a free non-symmetric operad.

A free non-symmetric operad describes type of algebras which have a set of generating operations and no relations between them. 
For instance, magmatic algebras are vector
spaces together with a bilinear product. There is a well known free non-symmetric operad, the Stasheff operad, built on Stasheff polytopes, see e.g. \cite{Stasheff63}. An algebra over the Stasheff operad is a vector space $V$ together with 
an $n$-linear product: $V^{\otimes n}\rto V$ for each $n$. From the point of view of homotopy theory, the category of operads is a Quillen category and free operads play  
an essential role in the homotopy category. One wants to replace an operad $\calP$
by a quasi-free resolution, that is, a morphism of operads $\mathcal{Q}\rightarrow\calP$ where $\mathcal Q$ is a free operad endowed with a differential realizing an isomorphism in homology. For instance, a quasi-free resolution of $\widetilde{\As}$, in the category of non-symmetric operads, is given by the Stasheff operad. Algebras over this operad 
are $A_\infty$-algebras (associative algebras up to homotopy). This gives us the motivation for studying whether a given symmetric operad is free as 
a non-symmetric operad or not.

\medskip

In this paper we prove that the operad pre-Lie is a free non-symmetric operad. Pre-Lie algebras are vector spaces together with a bilinear product satisfying the
relation $(x* y)* z-x*(y* z)=(x* z)* y-x*(z* y)$. The operad pre-Lie is based on labelled rooted trees which are of combinatorial interest. In the process of
proving the main result, we describe another operad denoted $\calT_\M$ also based on rooted trees and having the advantage of being the linearization of an operad in the category of sets. We prove that
it is a free non-symmetric operad. The link between the two operads is made via a gradation on labelled rooted trees.


\section{The pre-Lie operad and rooted trees}\label{S:LieT}

 We first recall the definition of the pre-Lie operad based
  on labelled rooted trees as in \cite{ChaLiv01}. For $n\in \mathbb N^*$,
the set $\{1,\ldots,n\}$ is denoted by $[n]$ and $[0]$ denotes the empty set.
The symmetric group on $k$ letters is denoted by $S_k$. There are many equivalent definitions of operads and we refer to
\cite{MSS02} for basics on operads.  We work over the ground field $\mathbf{k}$ and vector spaces are considered over $\mathbf k$. 
Here are the definitions needed for the sequel.

\begin{defi}\label{defop} A (reduced) {\sl non-symmetric} operad is a graded vector space $(\calP(n))_{n\geq 1}$,
with a unit $1\in\calP(1)=\mathbf k$,  together with
composition maps $\circ_i:\calP(n)\otimes\calP(m)\rto \calP(n+m-1)$ for $1\leq i\leq n$ satisfying the following relations: for $a\in\calP(n)$, $b\in\calP(m)$ and $c\in\calP(\ell)$
$$\begin{array}{lll}
(a\circ_i b)\circ_{j+i-1} c&=a\circ_i(b\circ_j c) ,& \hbox{for $1\le j\le m$,} \\
(a\circ _i b)\circ_j c&=(a\circ_{j}c)\circ_{i+\ell-1} b, &\textrm{for}\  j<i, \\
1\circ_1 a&=a, \\
a\circ_i 1&=a,
\end{array}$$
A non-trivial composition is a composition $a\circ_i b$ with $a\in\calP(n), b\in\calP(m)$ and $n,m>1$.

If in addition each $P(n)$ is acted on the right by the symmetric group $S_n$ and the composition maps are equivariant with respect to this action, then the collection
$(\calP(n))_n$ forms a {\sl symmetric operad}.
An {\sl algebra} over an operad $\calP$ is a vector space $X$ endowed with evaluation maps 
$$\begin{array}{cccc}
ev_n:&\calP(n)\otimes X^{\otimes n}&\rto& X\\
& p\otimes x_1\otimes\ldots\otimes x_n &\mapsto &p(x_1,\ldots,x_n)
\end{array}$$ 
compatible with the composition maps $\circ_i$:  for $p\in\calP(n),q\in\calP(m),x_i's\in X$ one has
$$(p\circ_i q)(x_1,\ldots,x_{n+m-1})=p(x_1\ldots,x_{i-1},q(x_i,\ldots,x_{i+m-1}),x_{i+m},\ldots,x_{n+m-1}).$$
If the operad is symmetric the evaluation maps are required to be equivariant with respect to the action of the symmetric group as follows:
$$(p\cdot \sigma)(x_1,\ldots,x_n)=p(x_{\sigma^{-1}(1)},\ldots,x_{\sigma^{-1}(n)}).$$
In the sequel an operad will always mean a reduced operad.

\end{defi}

\begin{defi} \label{def:trees}
In this paper we will consider two type of trees: planar rooted trees will represent 
the composition maps in a non-symmetric operad (see \ref{treeop}) and rooted trees will 
be the objects of our study (see \ref{S-rooted-tree}). Here are the definitions we will use in the sequel.

By a (planar) tree we mean a non empty finite connected contractible (planar) graph. All the trees considered
are rooted.

In the planar case some edges ({\sl external edges} or {\sl legs}) will have only one adjacent vertex; 
the other edges are called {\sl internal edges}. There is a distinguished leg called the {\sl root leg}. 
The other legs are called the leaves. The choice of a root induces a natural orientation of the 
graph from the leaves to the root. Any vertex has incoming edges and 
only one outgoing edge. The {\sl arity} of a vertex
is the number of incoming edges. A tree with no vertices of arity one is called {\sl reduced}. 
A planar rooted tree induces a structure of poset on the vertices, where $x<y$ if and only if there is 
an oriented path in the tree from $y$ to $x$. Let $x$ be a vertex of a planar rooted tree $T$. 
The {\sl full subtree $T^{(x)}$} of $T$ at $x$ is the subtree of $T$ containing all the vertices $y>x$ and
all their adjacent edges. The root leg of $T^{(x)}$ is the half edge with adjacent vertex $x$
induced by the unique outgoing edge of $x$. 
One represents a planar rooted tree like this:
$$ \raise -50pt
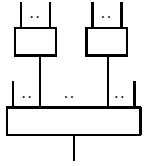
\hbox{ \begin{picture}(65,70)
       \put(20,6){\line(0,1){10}}   
       \put(-5,16){\framebox(50,10)[0,16]}
       \put(-2,46){\framebox(15,10)[0,16]}
       \put(25,46){\framebox(15,10)[0,16]}
       \put(-3,26){\line(0,1){10}}      \put(35,30){$\scriptscriptstyle ..$}       \put(43,26){\line(0,1){10}}      
       \put(0,56){\line(0,1){10}}    \put(3,60){$\scriptscriptstyle ..$}   \put(11,56){\line(0,1){10}}   
       \put(27,56){\line(0,1){10}}    \put(30,60){$\scriptscriptstyle ..$}   \put(38,56){\line(0,1){10}}   
        \put(16,30){$\scriptscriptstyle ..$}        \put(0,30){$\scriptscriptstyle ..$}
      \put(7,26){\line(0,1){20}}  
       \put(33,26){\line(0,1){20}}  
      \end{picture}}
$$

In the abstract case (non-planar trees) every edge is an internal edge. The {\sl root vertex} 
will be a distinguished vertex. The choice of a root induces a natural orientation of the 
graph towards the root. Any vertex has incoming edges and 
at most one outgoing edge. The other extremity of an incoming (outgoing) edge of the vertex $v$ is called an 
{\sl incoming (outgoing) vertex of the vertex $v$}.
The root vertex has no outgoing vertex.
A rooted tree induces a structure of poset on the vertices, where $x<y$ if and only if there is 
an oriented path in the tree from $y$ to $x$. A {\sl leave} is a maximal vertex for this order. 
The root is the only minimal vertex for this order.
Let $x$ be a vertex of a rooted tree $T$. 
The {\sl full subtree $T^{(x)}$} of $T$ derived from the vertex $x$ is the subtree of 
$T$ containing all the vertices $y>x$. The root of $T^{(x)}$ is $x$.
One represents a rooted tree like this:
$$\raise -10pt
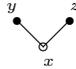
\hbox{ \begin{picture}(40,30)
       \put(10,6){\circle{3}} \put(0,16){\circle*{3}}  \put(20,16){\circle*{3}}
       \put(10,6){\line(-1,1){10}}  \put(10,6){\line(1,1){10}}
       \put(10,-1){$\scriptscriptstyle x$}   \put(-4,20){$\scriptscriptstyle y$} \put(20,20){$\scriptscriptstyle z$}  
      \end{picture}}$$ 

\end{defi}

\begin{rem} \label{treeop}{\sl Reduced planar tree of operations:} a convenient way to uniquely 
represent composition of operations in a non-symmetric operad $\calP$ is to use a planar rooted tree as in 
Definition~\ref{def:trees}. 
An element $a\in\calP(n)$ is represented by a planar rooted tree with a single vertex labelled by $a$
with $n$ incoming legs and 
a single outgoing leg:
$$ \raise -10pt
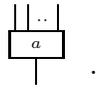
\hbox{ \begin{picture}(30,40)
       \put(10,6){\line(0,1){10}}   
       \put(0,16){\framebox(20,10)[0,16]{$\scriptscriptstyle a$}}
       \put(2,26){\line(0,1){10}}       \put(7,26){\line(0,1){10}}   \put(18,26){\line(0,1){10}}      
       \put(10,30){$\scriptscriptstyle ..$}
      \end{picture}} .
$$
The $n$ leaves are counted from left to right as $1$, $2$, \ldots, $n$. Now if we have $a\in\calP(n)$, $b\in\calP(m)$ and $1\le i\le n$ we 
represent the composition $a\circ_i b$ by the planar tree
$$ \raise -10pt
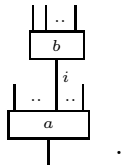
\hbox{ \begin{picture}(35,70)
       \put(10,6){\line(0,1){10}}   
       \put(-5,16){\framebox(30,10)[0,16]{$\scriptscriptstyle a$}}
       \put(3,46){\framebox(20,10)[0,16]{$\scriptscriptstyle b$}}
       \put(-3,26){\line(0,1){10}}       \put(13,26){\line(0,1){20}}   \put(23,26){\line(0,1){10}}      
       \put(16,30){$\scriptscriptstyle ..$}        \put(3,30){$\scriptscriptstyle ..$}
       \put(4,56){\line(0,1){10}}       \put(9,56){\line(0,1){10}}   \put(20,56){\line(0,1){10}}      
       \put(12,60){$\scriptscriptstyle ..$}
       \put(15,37){$\scriptscriptstyle i$}
      \end{picture}} .
$$
The resulting tree has $n+m-1$ leaves (counted from left to right) and represents an element of $\calP(n+m-1)$. 
The two first relations in Definition~\ref{defop} corresponds to the following two trees: for $a\in\calP(n)$, 
$b\in\calP(m)$ and $c\in\calP(\ell)$ we can have
$$ \raise -50pt
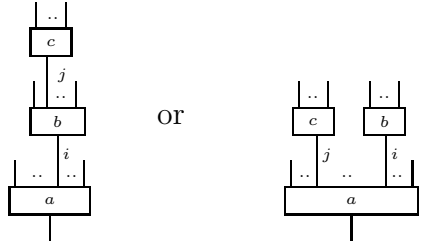
\hbox{ \begin{picture}(35,100)
       \put(10,6){\line(0,1){10}}   
       \put(-5,16){\framebox(30,10)[0,16]{$\scriptscriptstyle a$}}
       \put(3,46){\framebox(20,10)[0,16]{$\scriptscriptstyle b$}}
       \put(3,76){\framebox(15,10)[0,16]{$\scriptscriptstyle c$}}
       \put(-3,26){\line(0,1){10}}       \put(13,26){\line(0,1){20}}   \put(23,26){\line(0,1){10}}      
       \put(16,30){$\scriptscriptstyle ..$}        \put(3,30){$\scriptscriptstyle ..$}
       \put(4,56){\line(0,1){10}}       \put(9,56){\line(0,1){20}}   \put(20,56){\line(0,1){10}}      
       \put(5,86){\line(0,1){10}}        \put(16,86){\line(0,1){10}}      
       \put(12,60){$\scriptscriptstyle ..$}
       \put(15,37){$\scriptscriptstyle i$}
        \put(9,90){$\scriptscriptstyle ..$}
       \put(13,67){$\scriptscriptstyle j$}
       \end{picture}}  
 \quad\hbox{ or }\quad\qquad
       \raise -50pt\hbox{ \begin{picture}(65,70)
       \put(20,6){\line(0,1){10}}   
       \put(-5,16){\framebox(50,10)[0,16]{$\scriptscriptstyle a$}}
       \put(-2,46){\framebox(15,10)[0,16]{$\scriptscriptstyle c$}}
       \put(25,46){\framebox(15,10)[0,16]{$\scriptscriptstyle b$}}
       \put(-3,26){\line(0,1){10}}      \put(35,30){$\scriptscriptstyle ..$}       \put(43,26){\line(0,1){10}}      
       \put(0,56){\line(0,1){10}}    \put(3,60){$\scriptscriptstyle ..$}   \put(11,56){\line(0,1){10}}   
       \put(27,56){\line(0,1){10}}    \put(30,60){$\scriptscriptstyle ..$}   \put(38,56){\line(0,1){10}}   
        \put(16,30){$\scriptscriptstyle ..$}        \put(0,30){$\scriptscriptstyle ..$}
      \put(7,26){\line(0,1){20}}   \put(9,37){$\scriptscriptstyle j$}
       \put(33,26){\line(0,1){20}}   \put(35,37){$\scriptscriptstyle i$}
      \end{picture}}
$$
Each relation is obtained by writing down the two ways of interpreting the tree as a composition of operations. In general a planar tree $\TT(a_1,a_2,\ldots,a_k)$ with $k$ vertices labelled by elements
$a_i\in\calP(n_i)$ where $n_i$ is the number of incoming edges at the i$th$ vertex, corresponds to a unique composition of operations in $\calP$ 
independent of any relations.

The two last relations in Definition~\ref{defop} say that one can consider reduced trees (no vertices of arity 1)
for reduced operads to represent non-trivial composition maps.
  
Any full subtree of $\TT(a_1,a_2,\ldots,a_k)$ is completely determined 
by the position of its leaves; they form an interval $[p,q]$ where 
$1\le p\le q\le n_1+n_2+\cdots n_k-k+1$. A tree {\sl in position $[p,q]$} will mean the 
full subtree determined by the position $[p,q]$ of its leaves. If a full subtree in position $[p,q]$ has
a single vertex labelled by $a\in\calP(n)$ we identify this tree with the element $a\in\calP(n)$. It is clear that $n=q-p+1$.

Two trees of operations $\TT(a_1,a_2,\ldots,a_k)$ and $\YY(b_1,b_2,\ldots,b_s)$ are distinct if and only if
$\TT\not=\YY$ or there exists $i$ such that $a_i\not=b_i$.
\end{rem}

\smallskip

\begin{defi}\label{S-rooted-tree} Let $S$ be a set. An {\sl $S$-labelled rooted tree} is a non planar rooted 
tree as in Definition~\ref{def:trees} whose
  vertices are in bijection with $S$. If $S=[n]$, then we talk about $n$-labelled rooted trees and 
denote by $\calT(n)$ the set of those trees.  
It is acted on by the symmetric group by permuting the labels.

 The set $\calT(3)$ has for elements:
   \begin{equation}\label{e:T3} 
 \raise -10pt
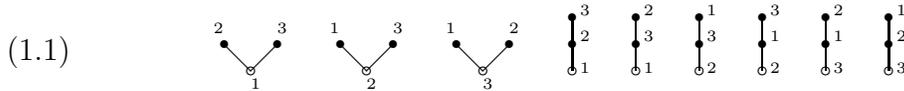
\hbox{ \begin{picture}(40,30)
       \put(10,6){\circle{3}}   \put(0,16){\circle*{3}}   \put(20,16){\circle*{3}}
       \put(10,6){\line(-1,1){10}}    \put(10,6){\line(1,1){10}}
       \put(10,-1){$\scriptscriptstyle 1$}   \put(-4,20){$\scriptscriptstyle 2$}   \put(20,20){$\scriptscriptstyle 3$}  
      \end{picture}} 
 \raise -10pt\hbox{  \begin{picture}(40,30)
       \put(10,6){\circle{3}} \put(0,16){\circle*{3}} \put(20,16){\circle*{3}}
       \put(10,6){\line(-1,1){10}}  \put(10,6){\line(1,1){10}}
       \put(10,-1){$\scriptscriptstyle 2$}   \put(-4,20){$\scriptscriptstyle 1$} \put(20,20){$\scriptscriptstyle 3$}  
      \end{picture}} 
 \raise -10pt\hbox{ \begin{picture}(40,30)
       \put(10,6){\circle{3}} \put(0,16){\circle*{3}}  \put(20,16){\circle*{3}}
       \put(10,6){\line(-1,1){10}}  \put(10,6){\line(1,1){10}}
       \put(10,-1){$\scriptscriptstyle 3$}   \put(-4,20){$\scriptscriptstyle 1$} \put(20,20){$\scriptscriptstyle 2$}  
      \end{picture}} 
 \raise -10pt\hbox{ \begin{picture}(20,30)
       \put(0,6){\circle{3}} \put(0,16){\circle*{3}} \put(0,26){\circle*{3}}
       \put(0,6){\line(0,1){20}}
       \put(3,5){$\scriptscriptstyle 1$} \put(3,17){$\scriptscriptstyle 2$} \put(3,27){$\scriptscriptstyle 3$}  
      \end{picture}} 
 \raise -10pt\hbox{ \begin{picture}(20,30)
       \put(0,6){\circle{3}}  \put(0,16){\circle*{3}} \put(0,26){\circle*{3}}
       \put(0,6){\line(0,1){20}}
       \put(3,5){$\scriptscriptstyle 1$}  \put(3,17){$\scriptscriptstyle 3$}  \put(3,27){$\scriptscriptstyle 2$}  
      \end{picture}} 
 \raise -10pt\hbox{ \begin{picture}(20,30)
       \put(0,6){\circle{3}}  \put(0,16){\circle*{3}} \put(0,26){\circle*{3}}
       \put(0,6){\line(0,1){20}} 
        \put(3,5){$\scriptscriptstyle 2$}  \put(3,17){$\scriptscriptstyle 3$} \put(3,27){$\scriptscriptstyle 1$}  
      \end{picture}} 
 \raise -10pt\hbox{ \begin{picture}(20,30)
       \put(0,6){\circle{3}}  \put(0,16){\circle*{3}} \put(0,26){\circle*{3}}
       \put(0,6){\line(0,1){20}}
       \put(3,5){$\scriptscriptstyle 2$}  \put(3,17){$\scriptscriptstyle 1$}  \put(3,27){$\scriptscriptstyle 3$}  
      \end{picture}} 
 \raise -10pt\hbox{ \begin{picture}(20,30)
       \put(0,6){\circle{3}} \put(0,16){\circle*{3}}  \put(0,26){\circle*{3}}
       \put(0,6){\line(0,1){20}}
       \put(3,5){$\scriptscriptstyle 3$}  \put(3,17){$\scriptscriptstyle 1$}  \put(3,27){$\scriptscriptstyle 2$}  
      \end{picture}} 
 \raise -10pt\hbox{ \begin{picture}(20,30)
       \put(0,6){\circle{3}}  \put(0,16){\circle*{3}} \put(0,26){\circle*{3}}
       \put(0,6){\line(0,1){20}}
       \put(3,5){$\scriptscriptstyle 3$}   \put(3,17){$\scriptscriptstyle 2$}   \put(3,27){$\scriptscriptstyle 1$}  
      \end{picture}} 
\end{equation}
In general $\calT(n)$ has $n^{n-1}$ elements (see \cite{BLL98} for more details). 
  \end{defi} We denote by ${\mathbf k}\calT(n)$ the $\mathbf k$-vector space spanned by $\calT(n)$.
  
\begin{theo}{\cite[theorem 1.9]{ChaLiv01}} The collection $({\mathbf k}\calT(n))_{n\geq 1}$ forms
  a symmetric operad, the operad pre-Lie denoted by $\PL$. Algebras over this operad are pre-Lie algebras, that is,
  vector spaces $L$ together with a product $*$ satisfying the
  relation
$$(x* y)* z-x*(y* z)=(x* z)* y-x*(z* y), \ \forall x,y,z\in L.$$
\end{theo}




\medskip

We recall the operad structure of $\PL$ as explained in \cite{ChaLiv01}.
A rooted tree is naturally oriented from the leaves to the root. The set $\In(T,i)$ of incoming vertices of a vertex $i$ is the set of all vertices $j$ such that $(j,i)$ is an edge oriented from $j$ to $i$.  
There is also at most one outgoing vertex of a vertex $i$, i.e. a vertex $r$ such that $(i,r)$ is an oriented edge from $i$ to $r$, depending whether $i$ is the root of $T$ or not. For  
$T\in\calT(n)$ and $S\in\calT(m)$, we define
$$T\circ_i S=\sum_{f:In(T,i)\rightarrow [m]} T\circ_i^f S,$$
where $ T\circ_i^f S$ is the rooted tree obtained by substituting the tree $S$ for the vertex $i$ in $T$. 
The outgoing vertex of $i$, if it exists, becomes the outgoing vertex of the root of $S$, whereas the 
incoming vertices of $i$ are grafted
on the vertices of $S$ according to the map $f$. The root of $T\circ_i^f S$ is the root of 
$T$ if $i$ is not the root of $T$, and it is the root of $S$ if $i$ is the 
root of $T$. There is also a relabelling of the vertices of $T$ and $S$ in $T\circ_i^f S$: we add $i-1$ to the labels of $S$ and $m-1$ to the ones of $T$ which are
greater than $i$.
Here is an example:

\begin{equation}\label{excompo}
\raise -10pt\hbox{ \begin{picture}(25,30)
       \put(0,16){\circle*{3}}    \put(20,16){\circle*{3}}     \put(10,6){\circle{3}}
       \put(0,16){\line(1,-1){10}}    \put(20,16){\line(-1,-1){10}}                                 
       \put(-3,20){$\scriptscriptstyle 1$}  \put(22,20){$\scriptscriptstyle 3$}  \put(10,-1){$\scriptscriptstyle 2$} 
       \end{picture}}\ \circ_2  
\raise -10pt\hbox{ \begin{picture}(25,30)
       \put(10,16){\circle*{3}}     \put(10,6){\circle{3}}
       \put(10,16){\line(0,-1){10}}    
       \put(10,20){$\scriptscriptstyle 1$}  \put(10,-1){$\scriptscriptstyle 2$} 
       \end{picture}} =
 \raise -10pt\hbox{ \begin{picture}(25,30)
       \put(0,16){\circle*{3}}    \put(20,16){\circle*{3}}     \put(10,6){\circle{3}}
       \put(0,16){\line(1,-1){10}}    \put(20,16){\line(-1,-1){10}}                                 
       \put(-3,20){$\scriptscriptstyle 1$}  \put(22,20){$\scriptscriptstyle 4$}  \put(10,-1){$\scriptscriptstyle 2$} 
       \end{picture}}\ \circ_2  
\raise -10pt\hbox{ \begin{picture}(25,30)
       \put(10,16){\circle*{3}}     \put(10,6){\circle{3}}
       \put(10,16){\line(0,-1){10}}    
       \put(10,20){$\scriptscriptstyle 2$}  \put(10,-1){$\scriptscriptstyle 3$} 
       \end{picture}} =
     \raise -10pt\hbox{ \begin{picture}(25,30)
       \put(10,16){\circle*{3}}   \put(0,26){\circle*{3}}   \put(20,26){\circle*{3}}      \put(10,6){\circle{3}}
       \put(10,16){\line(-1,1){10}}    \put(10,16){\line(1,1){10}}                           \put(10,6){\line(0,1){10}}        
       \put(12,13){$\scriptscriptstyle 2$}   \put(-2,30){$\scriptscriptstyle 1$}  
        \put(20,30){$\scriptscriptstyle 4$}                                           \put(10,-1){$\scriptscriptstyle 3$} 
       \end{picture}}\, 
       +\ \  \raise -10pt\hbox{ \begin{picture}(15,30)
       \put(0,16){\circle*{3}}   \put(0,26){\circle*{3}}   \put(10,16){\circle*{3}}      \put(0,6){\circle{3}}
       \put(0,16){\line(0,1){10}}    \put(0,6){\line(1,1){10}}                           \put(0,6){\line(0,1){10}}        
       \put(-8,13){$\scriptscriptstyle 2$}   \put(-2,30){$\scriptscriptstyle 1$}  
        \put(10,20){$\scriptscriptstyle 4$}                                           \put(0,-1){$\scriptscriptstyle 3$} 
       \end{picture}}\ 
 +\  \ \raise -10pt\hbox{ \begin{picture}(15,30)
       \put(0,16){\circle*{3}}   \put(0,26){\circle*{3}}   \put(10,16){\circle*{3}}      \put(0,6){\circle{3}}
       \put(0,16){\line(0,1){10}}    \put(0,6){\line(1,1){10}}                           \put(0,6){\line(0,1){10}}        
       \put(-8,13){$\scriptscriptstyle 2$}   \put(-2,30){$\scriptscriptstyle 4$}  
        \put(10,20){$\scriptscriptstyle 1$}                                           \put(0,-1){$\scriptscriptstyle 3$} 
       \end{picture}}\
  + \raise -10pt\hbox{ \begin{picture}(25,30)
       \put(0,16){\circle*{3}}   \put(10,16){\circle*{3}}   \put(20,16){\circle*{3}}      \put(10,6){\circle{3}}
       \put(0,16){\line(1,-1){10}}    \put(20,16){\line(-1,-1){10}}    \put(10,6){\line(0,1){10}}        
       \put(8,20){$\scriptscriptstyle 2$}   \put(-3,20){$\scriptscriptstyle 1$}  
        \put(21,20){$\scriptscriptstyle 4$}                                           \put(10,-1){$\scriptscriptstyle 3$} 
       \end{picture}}\,  
       \end{equation}

\bigskip


\section{A gradation on labelled rooted trees}\label{filtration}

We introduce a gradation on labelled rooted trees.  We prove that in the expansion 
of the composition  of two rooted trees in the operad pre-Lie  there is a unique rooted tree of maximal degree and a unique tree of minimal degree, yielding new non-symmetric operad structures on labelled rooted trees.

\begin{defi} Let $T$ be an $n$-labelled rooted tree. Let $\{a,b\}$ denote a pair of two adjacent vertices labelled
by $a$ and $b$. The degree of $\{a,b\}$ is $|a-b|$.
The {\sl degree of $T$} denoted by $\deg(T)$ is the sum of the degrees of its pairs of adjacent vertices.
For instance
$$ \deg(\raise -10pt\hbox{\begin{picture}(25,30)
       \put(0,16){\circle*{3}}    \put(20,16){\circle*{3}}     \put(10,6){\circle{3}}
       \put(0,16){\line(1,-1){10}}    \put(20,16){\line(-1,-1){10}}                                 
       \put(-3,20){$\scriptscriptstyle 1$}  \put(22,20){$\scriptscriptstyle 3$}  \put(10,-1){$\scriptscriptstyle 2$} 
       \end{picture}})=2 ,\quad 
 \deg(\raise -10pt\hbox{ \begin{picture}(25,30)
       \put(10,16){\circle*{3}}   \put(0,26){\circle*{3}}   \put(20,26){\circle*{3}}      \put(10,6){\circle{3}}
       \put(10,16){\line(-1,1){10}}    \put(10,16){\line(1,1){10}}                           \put(10,6){\line(0,1){10}}        
       \put(12,13){$\scriptscriptstyle 2$}   \put(-2,30){$\scriptscriptstyle 1$}  
        \put(20,30){$\scriptscriptstyle 4$}                                           \put(10,-1){$\scriptscriptstyle 3$} 
       \end{picture}})=4 ,\quad 
    \deg(\ \raise -10pt\hbox{ \begin{picture}(15,30)
       \put(0,16){\circle*{3}}   \put(0,26){\circle*{3}}   \put(10,16){\circle*{3}}      \put(0,6){\circle{3}}
       \put(0,16){\line(0,1){10}}    \put(0,6){\line(1,1){10}}                           \put(0,6){\line(0,1){10}}        
       \put(-8,13){$\scriptscriptstyle 2$}   \put(-2,30){$\scriptscriptstyle 4$}  
        \put(10,20){$\scriptscriptstyle 1$}                                           \put(0,-1){$\scriptscriptstyle 3$} 
       \end{picture}})=5,\quad
 \deg(\ \raise -10pt\hbox{ \begin{picture}(15,30)
       \put(0,16){\circle*{3}}   \put(0,26){\circle*{3}}   \put(10,16){\circle*{3}}      \put(0,6){\circle{3}}
       \put(0,16){\line(0,1){10}}    \put(0,6){\line(1,1){10}}                           \put(0,6){\line(0,1){10}}        
       \put(-8,13){$\scriptscriptstyle 2$}   \put(-2,30){$\scriptscriptstyle 1$}  
        \put(10,20){$\scriptscriptstyle 4$}                                           \put(0,-1){$\scriptscriptstyle 3$} 
       \end{picture}})=3
 $$
\end{defi}

\medskip

\begin{prop}\label{P-minmax} In the expansion of $T\circ_i S$ in the operad pre-Lie, there is a unique tree of minimal degree 
and a unique tree of maximal degree.
 \end{prop}

\medskip

For instance, in the equation (\ref{excompo}) the rooted tree of minimal degree 3 is $\ \ \raise -10pt\hbox{ \begin{picture}(15,30)
       \put(0,16){\circle*{3}}   \put(0,26){\circle*{3}}   \put(10,16){\circle*{3}}      \put(0,6){\circle{3}}
       \put(0,16){\line(0,1){10}}    \put(0,6){\line(1,1){10}}                           \put(0,6){\line(0,1){10}}        
       \put(-8,13){$\scriptscriptstyle 2$}   \put(-2,30){$\scriptscriptstyle 1$}  
        \put(10,20){$\scriptscriptstyle 4$}                                           \put(0,-1){$\scriptscriptstyle 3$} 
       \end{picture}}$
       and the one of maximal degree 5 is $\ \ \raise -10pt\hbox{ \begin{picture}(15,30)
       \put(0,16){\circle*{3}}   \put(0,26){\circle*{3}}   \put(10,16){\circle*{3}}      \put(0,6){\circle{3}}
       \put(0,16){\line(0,1){10}}    \put(0,6){\line(1,1){10}}                           \put(0,6){\line(0,1){10}}        
       \put(-8,13){$\scriptscriptstyle 2$}   \put(-2,30){$\scriptscriptstyle 4$}  
        \put(10,20){$\scriptscriptstyle 1$}                                           \put(0,-1){$\scriptscriptstyle 3$} 
       \end{picture}}$. The other ones are of degree 4.

\medskip

\noindent{\sl Proof--} Any tree in the expansion of $T\circ_i S$ writes $U_f:=T\circ_i^f S$ for some $f:\In(T,i)\rto[m]$. To compute the degree of $U_f$, we compute the degree of a pair of two adjacent vertices
$\{a,b\}$ in $U_f$. There are 4 cases to consider: a) the pair was previously in $S$ or
b) it was previously in $T$ and each vertex was different from $i$, or c) it was in $T$ of the form $\{i,j\}$ for $j\in\In(T,i)$ or d) if $i$ is not the root of $T$ it was of the form $\{i,k\}$ where $k$ is the outgoing vertex of $i$.

In case a) the degree of the pair in $U_f$ is the same as it was in $S$.

In case b), let $\{a',b'\}$ be the corresponding pair in $T$ before relabelling. The degree $d$ of the pair 
$\{a,b\}$ in $U_f$ is the same as the degree $d'$ of $\{a',b'\}$ except if $a'<i<b'$ or $b'<i<a'$,
where $d=d'+m-1$. Let $\gap(T,i)$ be the number of adjacent pairs of vertices in $T$ satisfying the latter condition.

In case c), let $\{i,j\}$ be the pair in $T$ which gives the pair $\{a,b\}$ in $U_f$. Let $d'$ be the degree of $\{i,j\}$. 
If $j<i$ then $\{a,b\}=\{f(j)+i-1,j\}$. Its degree $d$ is
minimal and equals $d'$ if $f(j)=1$. It is 
maximal and equals  $d'+m-1$ if $f(j)=m$.  
If $j>i$ then $\{a,b\}=\{f(j)+i-1,j+m-1\}$. Its degree $d$ is
minimal and equals $d'$ if $f(j)=m$. It is  
maximal and equals $d'+m-1$ if $f(j)=1$.

In case d), let $d'$ be the degree of $\{i,k\}$. 
If $k<i$ then
$\{a,b\}=\{s+i-1,k\}$ where $s$ is the label of the root of $S$. It has degree $d'+s-1$.  If 
$k>i$, then $\{a,b\}=\{s+i-1,k+m-1\}$ and has degree $(m-s)+d'$. 
Let $\epsilon(T,i,s)$ be $0,s-1,m-s$ according to the different situations, $0$ correponding to the one where $i$ is the root of $T$.

As a conclusion
\begin{multline}\label{E:degree}
\deg(T)+\deg(S)+\gap(T,i)(m-1)+\epsilon(T,i,s)\leq \deg(U_f)\leq \\
\deg(T)+\deg(S)+\gap(T,i)(m-1)+\epsilon(T,i,s)+|\In(T,i)|(m-1).
\end{multline}
There is a unique $f_\m$ such that $\deg(U_{f_\m})$ is minimal  and there is a unique $f_\M$ 
such that $\deg(U_{f_\M})$ is maximal:

\begin{equation}\label{fmin}
f_\m(k)=\begin{cases} 1 &\textrm{if\ } k<i, \\
 m &\textrm{if\ } k>i, \end{cases}
\end{equation}
\begin{equation}\label{fMax}
f_\M(k)=\begin{cases} m &\textrm{if\ } k<i, \\
 1 &\textrm{if\ } k>i, \end{cases}
\end{equation}
which ends the proof. \hfill $\Box$

\begin{theo}There are two different non-symmetric operad structures on the collection 
$({\mathbf k}\calT(n))_{n\geq 1}$ given by the composition maps
$T\circ_i^{f_\m} S$ on the one hand and $T\circ_i^{f_\M} S$ on the other hand where $f_\m$ and $f_\M$ were defined in
equations (\ref{fmin}) and (\ref{fMax}).
\end{theo}

\medskip

\noindent{\sl Proof--}  A  rooted tree $T$ is naturally oriented from its leaves to its root. Any edge is oriented and we denote by $(a,b)$ an edge oriented from the vertex $a$ to the vertex $b$. Let $E_T$ be the set of the oriented edges of the tree $T$.  For an integer $a\not=i$ we denote by
$\tilde a_i^m$ the integer $a$ if $a<i$ or $a+m-1$ if $a>i$. Given a map $f:\In(T,i)\rto [m]$, the set $E_{T\circ_i^f S}$ has different 
type of elements:  
\begin{itemize}
\item $(a+i-1,b+i-1)$ for $(a,b)\in E_S$; 
\item $(\tilde a_i^m,\tilde b_i^m)$ for $(a,b)\in E_T$ and $a,b\not=i$;
\item $(\tilde a_i^m,f(a)+i-1)$ for $(a,i)\in E_T$; 
\item $(i+s-1,\tilde b_i^m)$ for $(i,b)\in E_T$.
\end{itemize}

Let $T\in\calT(n), S\in\calT(m)$ and $U\in\calT(p)$. In order to avoid confusion, we denote by 
$f_\M^{i,p}$ the map sending $k<i$ to $p$ and $l>i$ to $1$.
We would like to compare the trees 
$$V_1=(T\circ_i^{f_\M^{i,m}} S)\circ_{j+i-1}^{f_\M^{j+i-1,p}} U\qquad \textrm{and}\qquad
V_2=T\circ_i^{f_\M^{i,m+p-1}}(S\circ_{j}^{f_\M^{j,p}} U):$$
\begin{itemize}
\item  In $V_1$ and $V_2$, any $(a,b)\in E_U$ converts to $(a+j+i-2,b+j+i-2)$.
\item In  $V_1$ and $V_2$, any $(a,b)\in E_S$ converts to $(\tilde a_j^p+i-1,\tilde b_j^p+i-1)$ if $a,b\not=j$, or converts to
$(\tilde a_j^p+i-1,f_\M^{j,p}(a)+i+j-2)$
if $b=j$ or converts to  $(j+i-1+u-1,\tilde b_j^p+i-1)$ if $a=j$.
\item In  $V_1$ and $V_2$, any $(a,b)\in E_T$ with $a,b\not=i$ converts to $(\tilde a_i^{p+m-1},\tilde b_i^{p+m-1})$. 
\item In  $V_1$ and $V_2$, any $(a,i)\in E_T$ converts to $(\tilde a_i^{p+m-1},f_\M^{i,m+p-1}(a)+i-1)$. 
\item In  $V_1$ and $V_2$, any $(i,b)\in E_T$ converts to $(i-1+\rac(S\circ_j U),\tilde b_i^{m+p-1})$, where $\rac(S\circ_j U)$ is the root
of $S\circ_j U$. More precisely 
$$\rac(S\circ_j U)=\begin{cases} s& \textrm{\ if } s<j\\ u+j-1 &\textrm{\  if\ } s=j \\ s+p-1 & \textrm{\ if\ } s>j.\end{cases}$$
\end{itemize}
The proof of 
$$(T\circ_i^{f_\M^{i,m}} S)\circ_{j}^{f_\M^{j,p}} U=(T\circ_j^{f_\M^{j,p}} U)\circ_{i+p-1}^{f_\M^{i+p-1,m}} S,\ \textrm{for}\ j<i$$
is similar and left to the reader. So is the proof with $f_\m$ instead of $f_\M$. \hfill$\Box$

 \medskip
 
The two operads on labelled rooted trees defined by the theorem are denoted by $\calT_\M$ and $\calT_\m$. Note that they are linearization of operads in the category of sets. Actually
the composition maps are defined at the level of the sets $\calT(n)$
and not only  at the level of the vector spaces ${\mathbf k}\calT(n)$. 
There is another operad built on rooted trees which has this property: the operad $\NAP$ encoding non-asociative permutative algebras in \cite{Livernet06}, in which $f_{\NAP}$ is the constant map with value the root of $S$. This operad has the advantage of being a symmetric operad.


\section{The operad pre-Lie is free as a non-symmetric operad}

We show that $\calT_\M$ is a free non-symmetric operad. Using Proposition~\ref{P-minmax}, we conclude that the operad pre-Lie is free as a non-symmetric operad. To this end we need to introduce some notation on rooted trees.


\begin{defi}
Given two ordered sets $S$ and $T$,
an order-preserving bijection  $\phi:S\rto T$ induces a natural 
bijection between the set of $S$-labelled rooted trees and the set of $T$-labelled rooted trees also denoted by $\phi$. 
A $T$-labelled rooted tree $X$ is {\sl isomorphic} to an $S$-labelled rooted tree $Y$ if $X=\phi(Y)$.
\end{defi}

Given a rooted tree $T\in\calT(n)$ and a subset $K\subseteq[n]$, we denote by  $T\big|_K$ the graph obtained from $T$ by keeping only the vertices of $T$ that are labelled by elements of $K$ and only the edges of $T$ that have two vertices labelled in $K$. Remark that each connected component of 
$T\big|_K$ is a rooted tree itself where the root is given by the unique vertex closest to the root of T in the component. Also, for $c\in[n]$ we denote by $T^{(c)}$ the full subtree of $T$ derived from the vertex labelled by $c$ (see Definition~\ref{def:trees}).
For example if $K=\{2,3,4,5,6\}\subset[7]$ and
$$T=\quad \raise -10pt
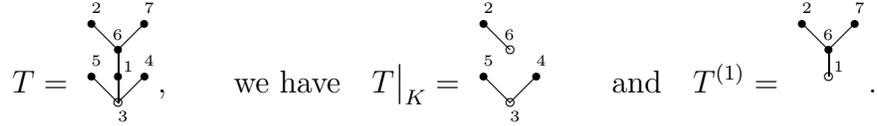
\hbox{ \begin{picture}(15,45)
       \put(0,16){\circle*{3}}   \put(0,26){\circle*{3}}   \put(10,16){\circle*{3}}      \put(0,6){\circle{3}} 
       \put(-10,16){\circle*{3}}   \put(-10,36){\circle*{3}}   \put(10,36){\circle*{3}}  
       \put(0,16){\line(0,1){10}}    \put(0,6){\line(1,1){10}}                           \put(0,6){\line(0,1){10}} 
          \put(0,6){\line(-1,1){10}}           \put(0,26){\line(1,1){10}}   \put(0,26){\line(-1,1){10}} 
       \put(2,18){$\scriptscriptstyle 1$}   \put(-2,30){$\scriptscriptstyle 6$}  
        \put(10,20){$\scriptscriptstyle 4$}    \put(0,-1){$\scriptscriptstyle 3$} 
        \put(-10,20){$\scriptscriptstyle 5$}    \put(10,40){$\scriptscriptstyle 7$} 
        \put(-10,40){$\scriptscriptstyle 2$} 
       \end{picture}},
       \hbox{\qquad we have\quad}
       T\big|_K=\quad \raise -10pt\hbox{ \begin{picture}(15,45)
      \put(0,26){\circle{3}}   \put(10,16){\circle*{3}}      \put(0,6){\circle{3}} 
       \put(-10,16){\circle*{3}}   \put(-10,36){\circle*{3}}  
       \put(0,6){\line(1,1){10}}             
          \put(0,6){\line(-1,1){10}}          \put(0,26){\line(-1,1){10}} 
  \put(-2,30){$\scriptscriptstyle 6$}  
        \put(10,20){$\scriptscriptstyle 4$}    \put(0,-1){$\scriptscriptstyle 3$} 
        \put(-10,20){$\scriptscriptstyle 5$} 
        \put(-10,40){$\scriptscriptstyle 2$} 
       \end{picture}}
              \hbox{\qquad and\quad}
       T^{(1)}=\quad \raise -10pt\hbox{ \begin{picture}(15,45)
       \put(0,16){\circle{3}}   \put(0,26){\circle*{3}}   
       \put(-10,36){\circle*{3}}   \put(10,36){\circle*{3}}  
       \put(0,16){\line(0,1){10}}                        
                  \put(0,26){\line(1,1){10}}   \put(0,26){\line(-1,1){10}} 
       \put(2,18){$\scriptscriptstyle 1$}   \put(-2,30){$\scriptscriptstyle 6$}  
      \put(10,40){$\scriptscriptstyle 7$} 
        \put(-10,40){$\scriptscriptstyle 2$} 
       \end{picture}}.
       $$

\medskip

For $1\le a<b\le n$, $T\in\calT_\M(n-b+a)$ and $S\in\calT_\M(b-a+1)$, let $X=T\circ_aS$.
Consider the interval $[a,b]=\{a,a+1,\ldots,b\}$,
clearly $X\big|_{[a,b]}$ is isomorphic to $S$ under the unique order-preserving bijection $[1,b-a+1]\to[a,b]$.
Let $a\le c\le b$ be the label of the root of  $X\big|_{[a,b]}$.
Remark that $X^{(c)}$ is obtained from $X\big|_{[a,b]}$ by grafting subtrees of $X$ at the vertices 
$a$ and $b$ only. 
We can then characterize trees $X$ that are obtained from a non-trivial composition $T\circ_aS$ as follows:

\begin{defi} \label{nprim}
A tree $X\in\calT_\M(n)$ is called {\sl decomposable} if there exists $1\le a<b\le n$  with $(a,b)\not=(1,n)$ such that
\begin{itemize}
\item[(i)] $X\big|_{[a,b]}$ is a rooted tree. Let $c$ be the label of its root.  One has $a\le c\le b$.
 \item[(ii)] One has $X^{(c)}|_{[a,b]}=X\big|_{[a,b]}$ and  $X^{(c)}$ is obtained from 
$X\big|_{[a,b]}$ by grafting subtrees of $X$ at the vertices 
$a$ and $b$ only. 
 \item[(iii)] All subtrees in $X^{(c)} - X\big|_{[a,b]}$ attached at $a$ have their root labelled in $[b+1,n]$.
  \item[(iv)] All subtrees in $X^{(c)} - X\big|_{[a,b]}$ attached at $b$ have their root labelled in $[1,a-1]$.
 \end{itemize}
 \end{defi}

It is clear from the discussion above and the definition of the operad $\calT_\M$ that $X$ is decomposable 
if and only if it is the result of a non-trivial composition. Consequently, we say that $X$ is {\sl indecomposable} 
if it is not decomposable. That is there is no $1\le a<b\le n$ such that (i)--(iv) are satisfied.
For example let
$$X=\quad \raise -10pt
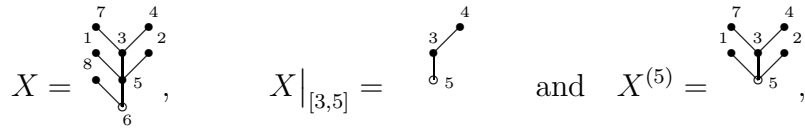
\hbox{ \begin{picture}(15,45)
       \put(0,16){\circle*{3}}   \put(0,26){\circle*{3}}   \put(10,26){\circle*{3}}      \put(0,6){\circle{3}} 
       \put(-10,16){\circle*{3}}   \put(-10,36){\circle*{3}}   \put(10,36){\circle*{3}}  
        \put(-10,26){\circle*{3}}  \put(0,16){\line(-1,1){10}}  \put(-15,30){$\scriptscriptstyle 1$}
       \put(0,16){\line(0,1){10}}    \put(0,16){\line(1,1){10}}                           \put(0,6){\line(0,1){10}} 
          \put(0,6){\line(-1,1){10}}           \put(0,26){\line(1,1){10}}   \put(0,26){\line(-1,1){10}} 
       \put(4,13){$\scriptscriptstyle 5$}   \put(-2,30){$\scriptscriptstyle 3$}  
        \put(13,30){$\scriptscriptstyle 2$}    \put(0,-1){$\scriptscriptstyle 6$} 
        \put(-15,20){$\scriptscriptstyle 8$}    \put(10,40){$\scriptscriptstyle 4$} 
        \put(-10,40){$\scriptscriptstyle 7$} 
       \end{picture}},
       \hbox{\qquad \quad}
       X\big|_{[3,5]}=\quad \raise -10pt\hbox{ \begin{picture}(15,45)
       \put(0,16){\circle{3}}   \put(0,26){\circle*{3}}     \put(10,36){\circle*{3}}  
       \put(0,16){\line(0,1){10}}      \put(0,26){\line(1,1){10}}  
       \put(4,13){$\scriptscriptstyle 5$}   \put(-2,30){$\scriptscriptstyle 3$}  
      \put(10,40){$\scriptscriptstyle 4$} 
       \end{picture}}
              \hbox{\qquad and\quad}
      X^{(5)}=\quad \raise -10pt\hbox{ \begin{picture}(15,45)
       \put(0,16){\circle{3}}   \put(0,26){\circle*{3}}   \put(10,26){\circle*{3}}    
      \put(-10,36){\circle*{3}}   \put(10,36){\circle*{3}}  
       \put(0,16){\line(0,1){10}}    \put(0,16){\line(1,1){10}}        
       \put(0,26){\line(1,1){10}}   \put(0,26){\line(-1,1){10}} 
        \put(-10,26){\circle*{3}}  \put(0,16){\line(-1,1){10}}  \put(-15,30){$\scriptscriptstyle 1$}
       \put(4,13){$\scriptscriptstyle 5$}   \put(-2,30){$\scriptscriptstyle 3$}  
        \put(13,30){$\scriptscriptstyle 2$}     \put(10,40){$\scriptscriptstyle 4$} 
        \put(-10,40){$\scriptscriptstyle 7$} 
       \end{picture}},
       $$
This tree $X$ is decomposable since for $1\le 3<5\le 8$ we have that $X\big|_{[3,5]}$ is a single tree and the subtrees of $X^{(5)}-X\big|_{[3,5]}$ are attached at $3$ and $5$ only. Moreover, the subtree attached  at 3 has root labelled by $7\in[6,8]$ and the subtrees attached  at 5 have roots labelled by $1,2\in[1,2]$. 
Indeed, in $\calT_\M$ we have
  $$ X = \quad \raise -10pt
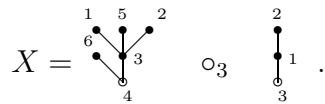
\hbox{ \begin{picture}(15,35)
       \put(0,16){\circle*{3}}   \put(0,26){\circle*{3}}   \put(10,26){\circle*{3}}      \put(0,6){\circle{3}} 
       \put(-10,16){\circle*{3}}    
        \put(-10,26){\circle*{3}}  \put(0,16){\line(-1,1){10}}  \put(-15,30){$\scriptscriptstyle 1$}
       \put(0,16){\line(0,1){10}}    \put(0,16){\line(1,1){10}}                           \put(0,6){\line(0,1){10}} 
          \put(0,6){\line(-1,1){10}}         
       \put(4,13){$\scriptscriptstyle 3$}   \put(-2,30){$\scriptscriptstyle 5$}  
        \put(13,30){$\scriptscriptstyle 2$}    \put(0,-1){$\scriptscriptstyle 4$} 
        \put(-15,20){$\scriptscriptstyle 6$}   
       \end{picture}} \quad \circ_3
       \quad \raise -10pt\hbox{ \begin{picture}(15,35)
       \put(0,16){\circle*{3}}   \put(0,26){\circle*{3}}    \put(0,6){\circle{3}} 
       \put(0,16){\line(0,1){10}}   \put(0,6){\line(0,1){10}} 
      \put(4,13){$\scriptscriptstyle 1$}   \put(-2,30){$\scriptscriptstyle 2$}  
      \put(0,-1){$\scriptscriptstyle 3$} 
       \end{picture}}.
       $$
The reader may check that the following are all the indecomposable trees of $\calT_\M$ up to arity 3:
  $$ \raise -10pt
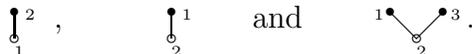
\hbox{ \begin{picture}(15,20)
       \put(0,16){\circle*{3}}   \put(0,6){\circle{3}}   \put(0,6){\line(0,1){10}} 
       \put(4,13){$\scriptscriptstyle 2$}   \  \put(0,-1){$\scriptscriptstyle 1$} 
       \end{picture}}, \quad\qquad
       \raise -10pt\hbox{ \begin{picture}(15,20)
       \put(0,16){\circle*{3}}   \put(0,6){\circle{3}}   \put(0,6){\line(0,1){10}} 
       \put(4,13){$\scriptscriptstyle 1$}   \  \put(0,-1){$\scriptscriptstyle 2$} 
       \end{picture}} \quad \hbox{ and }\qquad
       \quad \raise -10pt\hbox{ \begin{picture}(15,20)
       \put(10,16){\circle*{3}}   \put(-10,16){\circle*{3}}    \put(0,6){\circle{3}} 
       \put(0,6){\line(-1,1){10}}   \put(0,6){\line(1,1){10}} 
      \put(-16,13){$\scriptscriptstyle 1$}   \put(13,13){$\scriptscriptstyle 3$}  
      \put(0,-1){$\scriptscriptstyle 2$} 
       \end{picture}}\ .
       $$

\begin{theo}\label{T:Tmax}
The non-symmetric operad  $\calT_\M$ is a free non-symmetric operad.
\end{theo}


\begin{proof} If $\calT_\M$ is not free, then for some $n$ there is a tree $X\in\calT_\M(n)$ with two distinct constructions from indecomposables. In Remark~\ref{treeop}, a non-trivial composition of 
operations is completely determined 
by a unique reduced planar rooted tree. We then have that $X=\TT(T_1,T_2,\ldots,T_r)=\YY(S_1,S_2,\ldots,S_k)$ where $ T_1,\ldots,T_r,S_1,\ldots,S_k$ are indecomposables and $\TT(T_1,T_2,\ldots,T_r)$ and $\YY(S_1,S_2,\ldots,S_k)$ are two distinct trees of operations in $\calT_\M$ with $r , k>1$.

The tree $X=\TT(T_1,T_2,\ldots,T_r)$ is decomposable (by assumption $r\geq 2$). 
We can find $1\le a<b\le n$, such that $X\big|_{[a,b]}$ is isomorphic to a single $T_i$ in  position 
$[a,b]$ in $\TT(T_1,T_2,\ldots,T_r)$. Moreover $X\big|_{[a,b]}$ satisfies (i)--(iv) of 
Definition~\ref{nprim}. 

If $X\big|_{[a,b]}$ is also isomorphic to a tree $S_j$ 
in position $[a,b]$ in $\YY(S_1,S_2,\ldots,S_k)$, 
then we replace $X$ by the smaller tree in $\calT_\M(n-b+a)$ that we obtain by removing $T_i$ in $\TT(T_1,T_2,\ldots,T_r)$ and removing $S_j$ in $\YY(S_1,S_2,\ldots,S_k)$. Clearly, this new smaller $X$ has two distinct constructions from indecomposables. We can thus assume that $X\big|_{[a,b]}$ is not isomorphic to a single $S_j$ in position $[a,b]$ in $\YY(S_1,S_2,\ldots,S_k)$. 

We now study how $X\big|_{[a,b]}$ overlaps in the position $[a,b]$ of $\YY(S_1,S_2,\ldots,S_k)$. Remark first that since all $S_j$ are indecomposables, the interval $[a,b]$ cannot be part of a single $S_j$ of $\YY(S_1,S_2,\ldots,S_k)$. Indeed, that would imply that $S_j$ would contain a subtree satisfying Definition~\ref{nprim} which would be a contradiction.

We may assume that $a>1$. To see this, assume that the only sub-interval $[a,b]\subset[1,n]$ such that $X\big|_{[a,b]}$ is isomorphic to a single $T_i$ in  position 
$[a,b]$ in $\TT(T_1,T_2,\ldots,T_r)$ is such that $a=1$. Assume moreover that the only sub-interval $[a',b']\subset[1,n]$ such that $X\big|_{[a',b']}$ is isomorphic to a single $S_j$ in  position 
$[a',b']$ in $\YY(S_1,S_2,\ldots,S_k)$ is such that $a'=1$. Since $S_j$ is indecomposable, we must have $b> b'$. Similarly, since $T_i$ is indecomposable, we must have $b<b'$. This implies that $b=b'$ and $T_i=S_j$. This possibility was excluded above. So we must have $a>1$ or $a'>1$. In the case where $a=1$ and $a'>1$  we could just interchange the role of $\TT(T_1,T_2,\ldots,T_r)$ and $\YY(S_1,S_2,\ldots,S_k)$ and assume that we have $a>1$.

Now, since $T_i$ is indecomposable, there is no subinterval $[c,d]\subseteq [a,b]$ such that 
$X\big|_{[c,d]}$ is isomorphic to a full subtree of operations $\YY'(S_{j_1},S_{j_2},\ldots,S_{j_\ell})$.
Assume we can find $c<a\le d<b$ such that $X\big|_{[c,d]}\cong \YY'(S_{j_1},S_{j_2},\ldots,S_{j_\ell})$ 
satisfies the Definition~\ref{nprim}.

The graph $X\big|_{[a,d]}$ is contained in the trees $X\big|_{[a,b]}$ and $X\big|_{[c,d]}$.
Let $e$ be the label of the root of $X\big|_{[a,b]}$ and $f$ be the label of the root of $X\big|_{[c,d]}$. 
The two full subtrees $X^{(e)}$ and $X^{(f)}$ both contain  $X\big|_{[a,d]}$. This implies that either $X^{(f)}$ 
is fully contained in  $X^{(e)}$, or $X^{(e)}$ is fully contained in  $X^{(f)}$. 

Let us assume that $X^{(f)}$ is fully contained in  $X^{(e)}$, that means $X\big|_{[a,b]}$ and $X\big|_{[c,d]}$ 
are both subtrees of $X^{(e)}$. From Definition~\ref{nprim}, we know that
$X^{(e)}$ is obtained from $X\big|_{[a,b]}$ by graphting subtrees of $X$ at the vertices $a$ and $b$ only. The vertex $c$ is in $X^{(e)}$ but not in $X\big|_{[a,b]}$. It is part of a subtree attached to $a$ or $b$. Since $c$
is part of a subtree with root $f$ one has $f\not\in]a,b[$. The vertex $f$ is
$a$ (can not be $b$ since $f\le d$) or is attached 
to $a$ or $b$. If $f$ is attached to $b$ then there is a path $c\rightarrow f \rightarrow b$. 
The tree $X|_{[c,d]}$
has its root labelled by $f$ so there is a path $d\rightarrow f$. The tree $X|_{[a,b]}$ contains the vertices 
$b$ and $d$ and any path from $d$ to $b$ so  there is a path $d\rightarrow f\rightarrow b$ in $X|_{[a,b]}$.
Hence $f=a$ for $f\not\in\; ]a,b]$.
    As a conclusion $c$ is part of a subtree attached to $a$. By (iii) of 
Definition~\ref{nprim} applied to the tree $X\big|_{[a,b]}$, the subtree
must have a root $r\in[b+1,n]$. This is a contradiction, the root  $r$ is part of 
any path joining $a$ and $c$ and $r\not\in[c,d]$, hence not in $X\big|_{[c,d]}$. The case where $X^{(e)}$ is 
fully contained in  $X^{(f)}$ is argued similarly, using condition (iv) of Definition~\ref{nprim}, and leads 
to a contradiction as well. 

The same argument holds in case we can find $a<c\le b<d$.

The only case remaining is that
the interval $[p,q]$ associated to any full subtree $\YY(S_1,\ldots,S_k)^{(S_j)}$  
of $\YY(S_1,\ldots,S_k)$, satisfies $[a,b]\cap [p,q]=\emptyset$ or $[a,b]\subset [p,q]$. 
There is at least one interval satisfying  $[a,b]\subset [p,q]$ (take the full tree  $\YY(S_1,\ldots,S_k)$ and 
$[p,q]=[1,n]$). Let $[p,q]$ be the smallest interval such that $[a,b]\subset [p,q]$ and let
$\YY(S_1,\ldots,S_k)^{(S_{j})}=\YY'(S_{i_1},\ldots,S_{i_l})$ be the full subtree it determines. Its root is labelled by 
$S_j$. The interval $[u,v]$ associated to any proper full subtree of $\YY'(S_{i_1},\ldots,S_{i_l})$ 
satisfies
 $[a,b]\cap [u,v]=\emptyset$. Consequently $X|_{[a,b]}$ is isomorphic to $S_j|_{[\alpha,\beta]}$ for some interval
$[\alpha,\beta]$ isomorphic to $[a,b]$. This is impossible since $X$ satisfies the conditions of 
Definition~\ref{nprim} and $S_j$ is indecomposable.

We must conclude that  $\calT_\M$ is free. \end{proof}

\medskip

\begin{rem} The non-symmetric operads $\calT_\m$ and $\NAP$ are not free. 
Indeed, in the operad $\calT_\m$ one has the following relation:
$$ \raise -10pt\hbox{ \begin{picture}(25,30)
       \put(10,16){\circle*{3}}     \put(10,6){\circle{3}}
       \put(10,16){\line(0,-1){10}}    
       \put(10,20){$\scriptscriptstyle 2$}  \put(10,-1){$\scriptscriptstyle 1$} 
       \end{picture}}\ \circ_1  
\raise -10pt\hbox{ \begin{picture}(25,30)
       \put(10,16){\circle*{3}}     \put(10,6){\circle{3}}
       \put(10,16){\line(0,-1){10}}    
       \put(10,20){$\scriptscriptstyle 2$}  \put(10,-1){$\scriptscriptstyle 1$} 
       \end{picture}} =
      \raise -10pt\hbox{ \begin{picture}(25,30)
       \put(10,16){\circle*{3}}     \put(10,6){\circle{3}}
       \put(10,16){\line(0,-1){10}}    
       \put(10,20){$\scriptscriptstyle 2$}  \put(10,-1){$\scriptscriptstyle 1$} 
       \end{picture}}\ \circ_2  
\raise -10pt\hbox{ \begin{picture}(25,30)
       \put(10,16){\circle*{3}}     \put(10,6){\circle{3}}
       \put(10,16){\line(0,-1){10}}    
       \put(10,20){$\scriptscriptstyle 2$}  \put(10,-1){$\scriptscriptstyle 1$} 
       \end{picture}} =\ \   \raise -10pt\hbox{ \begin{picture}(15,30)
       \put(0,16){\circle*{3}}   \put(0,26){\circle*{3}}    \put(0,6){\circle{3}}
       \put(0,16){\line(0,1){10}}   \put(0,6){\line(0,1){10}}        
       \put(-8,13){$\scriptscriptstyle 2$}   \put(-2,30){$\scriptscriptstyle 3$}  
         \put(0,-1){$\scriptscriptstyle 1$} 
       \end{picture}}$$

And in the operad $\NAP$ one has the following relation

$$ \raise -10pt\hbox{ \begin{picture}(25,30)
       \put(10,16){\circle*{3}}     \put(10,6){\circle{3}}
       \put(10,16){\line(0,-1){10}}    
       \put(10,20){$\scriptscriptstyle 2$}  \put(10,-1){$\scriptscriptstyle 1$} 
       \end{picture}}\ \circ_1  
\raise -10pt\hbox{ \begin{picture}(25,30)
       \put(10,16){\circle*{3}}     \put(10,6){\circle{3}}
       \put(10,16){\line(0,-1){10}}    
       \put(10,20){$\scriptscriptstyle 1$}  \put(10,-1){$\scriptscriptstyle 2$} 
       \end{picture}} =
      \raise -10pt\hbox{ \begin{picture}(25,30)
       \put(10,16){\circle*{3}}     \put(10,6){\circle{3}}
       \put(10,16){\line(0,-1){10}}    
       \put(10,20){$\scriptscriptstyle 1$}  \put(10,-1){$\scriptscriptstyle 2$} 
       \end{picture}}\ \circ_2  
\raise -10pt\hbox{ \begin{picture}(25,30)
       \put(10,16){\circle*{3}}     \put(10,6){\circle{3}}
       \put(10,16){\line(0,-1){10}}    
       \put(10,20){$\scriptscriptstyle 2$}  \put(10,-1){$\scriptscriptstyle 1$} 
       \end{picture}} =
     \raise -10pt\hbox{ \begin{picture}(25,30)
       \put(0,16){\circle*{3}}    \put(20,16){\circle*{3}}     \put(10,6){\circle{3}}
       \put(0,16){\line(1,-1){10}}    \put(20,16){\line(-1,-1){10}}                                 
       \put(-3,20){$\scriptscriptstyle 1$}  \put(22,20){$\scriptscriptstyle 3$}  \put(10,-1){$\scriptscriptstyle 2$} 
       \end{picture}}$$

\end{rem}

\begin{rem}\label{hilb} Let ${\mathbf k}\calT_\M^0(n)$ denote the $\mathbf k$-vector space spanned by the indecomposables of $\calT_\M(n)$ ($n>1)$ and let $\beta_n$ be its dimension. 
Let $\alpha(x)=\sum_{n\geq 1}\alpha_n x^n$ be the Hilbert series associated to the free non-symmetric operad generated by the vector spaces  ${\mathbf k}\calT_\M^0(n)$.
It is well known (see e.g. \cite{SalTau08}) that one has the identity
$$\beta(\alpha(x))+x=\alpha(x),$$
where $\beta(x)=\sum_{n\geq 2}\beta_nx^n$. Theorem \ref{T:Tmax} implies that $\alpha_n=n^{n-1}$. As a consequence,
we get that the Hilbert series for indecomposable of $\calT_\M$ is
  $${\mathcal H}_{\calT_\M^0}(x)=\sum_{n\ge 2} \dim \big({\mathbf k}\calT_\M^0(n)\big)\, x^n
     = 2x^2+ x^3+ 14x^4+146x^5+$$
     $$\ \qquad\qquad +1994x^6+ 32853x^7+ 630320x^8+13759430x^9+\cdots.
$$
\end{rem}

\begin{coro}The non-symmetric operad pre-Lie is a free non-symmetric operad. 
\end{coro}

\begin{proof} Let $\mathcal F$ be the free non-symmetric operad on indecomposable trees. 
By the universal property of  $\mathcal F$, there is a unique morphism of operads
$$\phi:\mathcal F\rightarrow \PL$$
extending the inclusion of indecomposable trees in $\PL$.
We prove that this map is surjective by induction on the degree of a tree.
Trees of degree 1 are indecomposables (see Definition~\ref{nprim}). Let $t\in\PL(n)$ be a tree of degree $k\geq n-1$.
If $t$ is indecomposable then $t=\phi(t)$. If $t$ is decomposable there are trees $u\in\PL(r), v\in\PL(s)$, with
$r,s<n$ such that $t=u\circ^{f_\M}_i v$ in $\calT_\M$. By Proposition~\ref{P-minmax} one has in $\PL$
$$u\circ_i v= t+\sum_j t_j$$ 
where $t_j\in\PL(n)$ has degree $k_j<k$. From equation (\ref{E:degree}) we deduce that
the degrees of $u$ and $v$ are also lower than $k$. By induction, the trees $u,v$ and $t_j's$ are in the image of 
$\phi$, so is $t$. Thus, the operad morphism $\phi$ is surjective.
Theorem~\ref{T:Tmax} implies that the vector spaces $\mathcal F(n)$ and 
$\PL(n)$ have the same dimension, thus the operad morphism $\phi$ is an isomorphism.\end{proof}


\begin{rem}
The Hilbert Series for the free non-symmetric operad on indecomposables and the operad $\PL$ are the same as in Remark~\ref{hilb}.
\end{rem}


\bibliographystyle{plain} 
\bibliography{bibliojuin08.bib}

\def\cprime{$'$} \def\cprime{$'$} \def\cprime{$'$} \def\cprime{$'$}
\begin{thebibliography}{1}

\bibitem{BLL98}
Fran\c{c}ois Bergeron, Gilbert Labelle, and Pierre Leroux.
\newblock {\em Combinatorial species and tree-like structures}, volume~67 of
  {\em Encyclopedia of Mathematics and its Applications}.
\newblock Cambridge University Press, Cambridge, 1998.
\newblock Translated from the 1994 French original by Margaret Readdy, With a
  foreword by Gian-Carlo Rota.

\bibitem{ChaLiv01}
Fr{\'e}d{\'e}ric Chapoton and Muriel Livernet.
\newblock Pre-{L}ie algebras and the rooted trees operad.
\newblock {\em Internat. Math. Res. Notices}, 8:395--408, 2001.

\bibitem{Livernet06}
Muriel Livernet.
\newblock A rigidity theorem for pre-{L}ie algebras.
\newblock {\em J. Pure Appl. Algebra}, 207(1):1--18, 2006.

\bibitem{MSS02}
Martin Markl, Steve Shnider, and Jim Stasheff.
\newblock {\em Operads in algebra, topology and physics}, volume~96 of {\em
  Mathematical Surveys and Monographs}.
\newblock American Mathematical Society, Providence, RI, 2002.

\bibitem{SalTau08}
Paolo Salvatore and Roberto Tauraso.
\newblock The operad {L}ie is free.
\newblock arXiv:0802.3010, 2008.

\bibitem{Stasheff63}
James~Dillon Stasheff.
\newblock Homotopy associativity of {$H$}-spaces. {I}, {II}.
\newblock {\em Trans. Amer. Math. Soc. 108 (1963), 275-292; ibid.},
  108:293--312, 1963.

\end{thebibliography}

\bigskip
\end{document}